\documentclass{article}

\usepackage{amsfonts,amssymb}

\begin{document}



\author{Ya.\,V.~Bazaikin, O.\,A.~Bogoyavlenskaya}

\title{Complete Riemannian $G_2$ Holonomy Metrics on Deformations of Cones over $S^3\times S^3$}


\maketitle

\begin{abstract}
Complete Riemannian metrics with holonomy group $G_2$ are constructed on the manifolds obtained by deformations of cones over $S^3 \times S^3$.

\end{abstract}


\section{Introduction}

This article is a sequel to the works \cite{Baz,Baz2,Baz3,Malk} and is dedicated to studying the Riemannian manifolds with holonomy  group $G_2$. Recently, this problem has attracted considerable interest; moreover, though the most important problem is to study the compact manifolds admitting such metrics, still the studying of non-compact manifolds  (for the most part,  the vector bundle spaces) with complete Riemannian $G_2$-holonomy  metrics is quite logical. This is explained by the fact that in the latter case one can succeed, as a rule, in setting a $G_2$-structure explicitly and writing the equations which guarantee its being parallel.  In addition, if the symmetry group of the considered $G_2$-structure is large enough, then the problem is reduced to a system of ordinary differential equations, which allows either to find explicit solutions (in contrast to the compact case), or to study them qualitatively.
The main idea of the article has been already used in \cite{Baz,Baz2,Baz3} for constructing  complete Riemannian  metrics with holonomy group $Spin(7)$; it consists in the following:  the standard conic metric is considered over a Riemannian manifold with a special geometry. Then the deformation of this metric depends on a certain number of functional parameters which allow defining explicitly a      
$G_2$ (or $Spin(7)$) structure. In the present work we propose (following  \cite{Brandhuber})  to consider the space $M=S^3 \times S^3$ as such base of the cone. Then the conic metric can be written as 
$$
d\bar{s}^2=dt^2+\sum_{i=1}^3 A_i(t)^2 \left( \eta_i+\tilde{\eta_i}\right)^2 
+ \sum_{i=1}^3 B_i(t)^2 \left( \eta_i-\tilde{\eta_i}\right)^2,
$$
where  $\eta_i, \tilde{\eta_i}$  is the standard  coframe of $1$-forms,  whereas the functions  $A_i(t), B_i(t)$  define a deformation of  the cone singularity.  In the paper \cite{Brandhuber}  a system of differential equations is written down, which guarantees that the metric $d\bar{s}^2$ has the holonomy group containing in $G_2$. In \cite{Brandhuber} a particular solution of this system is found, which corresponds to a metric with the holonomy group $G_2$ on  $S^3 \times \mathbb{R}^4$.  Let us indicate that in the papers    \cite{CGLP,CGLP2,B,CGLP3,CCGLPW}   more general metrics on the cones over $S^3\times S^3$ have  also been studied; however, in the situation considered by us no other examples have been discovered besides  the example from  \cite{Brandhuber} and the classical example from \cite{Br-Sal}.
In the proposed work we continue to study this class of metrics, while setting   
$A_2=A_3$, $B_2=B_3$  and considering boundary conditions different from that from  
\cite{Brandhuber}. This yields the metrics with a different topological structure.  Namely, we require 
that at the vertex of the cone  only the function $B_1$  turns to zero. This results in that the Riemannian 
metric $d\bar{s}^2$  is defined on $H^4 \times S^3$, where $H$ is the space of canonical complex linear 
bundle over $S^2$, where $H^4$ is its fourth tensor power. Note that in  \cite{CGLP3}  numerical 
investigation was conducted (using the development of the solution of the basic system into the Taylor 
series), which yielded some arguments in favor of existence of the metrics constructed by us.
The main result of the paper is formulated in the following theorem:

\vskip0.2cm

{\bf Theorem.} {\it There exists a one-parametric family of mutually non-homothetic complete 
Riemannian metrics of the form $d\bar{s}^2$  with holonomy group  $G_2$ on $H^4 \times S^3$,
whereas the metrics can be parameterized by the set of  initial data  $(A_1(0)$, $A_2(0)$, $B_1(0)$, 
$B_2(0))=(\mu, \lambda,0,\lambda)$, where  $\lambda, \mu>0$  and   
$\mu^2+\lambda^2=1$.

For $t \rightarrow \infty$  the metrics of this family are approximated arbitrary closely by the direct  
product  $S^1 \times C(S^2 \times S^3)$,  where $C(S^2 \times S^3)$  is the cone over the product of 
spheres. Moreover, the sphere $S^2$ arises as factorization of the diagonally embedded in $S^3 \times 
S^3$ three-dimensional sphere with respect to the circle action corresponding to the vector field 
$\xi^1+\tilde{\xi}_1$
}

\section{$G_2$-structure on the cone over $S^3 \times S^3$}

Consider the Lie group $G=SU(2)$ with the standard bi-invariant metric 
$$
\langle X,Y \rangle = -\mbox{tr } (X Y),
$$
where  $X,Y \in {\bf su}(2)$.
Let us consider three Killing vector fields on $G$:
$$
\xi^1=\left( \begin{array}{cc} i & 0 \\ 0 & -i \\ \end{array} \right), \
\xi^2=\left( \begin{array}{cc} 0 & 1 \\ -1 & 0 \\ \end{array} \right), \
\xi^3=\left( \begin{array}{cc} 0 & i \\ i & 0 \\ \end{array} \right).
$$
It is not difficult to see that they satisfy the relations 
$$
[\xi^i,\xi^{i+1}]=2 \xi^{i+2},
$$
where the indices $i=1, 2, 3$ are reduced modulo $3$. Let  $\eta_1, \eta_2, 
\eta_3$  be the dual basis of 1-forms, that is, 
 $\eta_i(\xi^j)=\delta_i^j$. 
Then 
$$
d\eta_i=- 2 \eta_{i+1} \wedge \eta_{i+2}.
$$
Let  $M=G \times G$, then on $M$ there arise  $6$ Killing fields  
$\xi^i$, $\tilde{\xi}^i$, $i=1, 2, 3$, which are tangent to the first and second factor, respectively, and 
$6$  dual 1-forms  $\eta_i$, $\tilde{\eta}_i$.
Consider the cone  $\overline{M}=\mathbb{R}_+ \times M$  with the metric 
$$
d\bar{s}^2=dt^2+\sum_{i=1}^3 A_i(t)^2 \left( \eta_i+\tilde{\eta}_i\right)^2 
+ \sum_{i=1}^3 B_i(t)^2 \left( \eta_i-\tilde{\eta}_i\right)^2,
$$
where $A_i(t)$  and $B_i(t)$ are some positive functions  defining a deformation of the standard 
conic metric.

Introducing the orthonormal coframe 
$$
\begin{array}{ll}
e^1=A_1 \left( \eta_1+\tilde{\eta}_1 \right), & e^4= B_1 \left( 
\eta_1-\tilde{\eta}_1 \right), \\
e^2=A_2 \left( \eta_2+\tilde{\eta}_2 \right), & e^5= B_2 \left( 
\eta_2-\tilde{\eta}_2 \right), \\
e^3=A_3 \left( \eta_3+\tilde{\eta}_3 \right), & e^6= B_3 \left( 
\eta_3-\tilde{\eta}_3 \right), \\
e^7=dt .& \\
\end{array}
$$
we define the following  3-form:
$$
\Psi=e^{564}+e^{527}+e^{513}+e^{621}+e^{637}+e^{432}+e^{417},
$$
where  $e^{ijk}=e^i \wedge e^j \wedge e^k$.  The form  $\Psi$ defines a  
$G_2$-structure on  $\overline{M}$, which is parallel  provided the following equations hold:
$$
d\Psi=0, d*\Psi=0.\eqno{(1)}
$$
In the present work we consider a particular case when 
$A_2=A_3$, 
$B_2=B_3$.

\vskip0.2cm

{\bf Lemma 1. } {\it  Equations  (1) are equivalent to the following system of ordinary differential 
equations:

$$
\begin{array}{l}
\frac{dA_1}{dt}=\frac{1}{2}\left(\frac{A_1^2}{A_2^2}-\frac{A_1^2}{B_2^2}\right)\\
\frac{dA_2}{dt}=\frac{1}{2}\left(\frac{B_2^2-A_2^2+B_1^2}{B_1B_2}-\frac{A_1}{A_2}\right)\\
\frac{dB_1}{dt}=\frac{A_2^2+B_2^2-B_1^2}{A_2B_2}\\
\frac{dB_2}{dt}=\frac{1}{2}\left(\frac{A_2^2-B_2^2+B_1^2}{A_2B_1}+\frac{A_1}{B_2}\right)
\end{array}\eqno(2)
$$

}

\vskip0.2cm

For $t=0$ we have a conic singularity of the space $\overline{M}$  which can be resolved by setting the initial values of the functions  $A_i, B_i$. At that, there appear, up to symmetry of system (2), two 
types of singularity resolution listed below (compare with \cite{Baz, Baz2, Baz3}).

Type 1. $A_i(0)=0, B(0)\neq 0$. In this case, a collapse takes place of the integral  three-dimensional 
spheres generated by the vector fields $\xi^i+\tilde{\xi}^i$.  These spheres are the orbits of a free action 
of $G$ on $M$, defined by the relation $h \in G: (g_1,g_2) \mapsto (hg_1, hg_2)$. 
It can be demonstrated that in this case the metric $d\bar{s}^2$ on $\overline{M}$ 
can be continued to a space ${\cal M}$, homeomorphic to  $S^3 \times 
\mathbb{R}^4$. Here we omit the details because this case is not investigated in the present article.

Type 2. $B_1(0)=0$, $B_2(0)\neq 0, A_i(0) \neq 0$. Consider a free action of the group $U(1)=S^1$ on $M$:
$$
z \in U(1): (U,V) \mapsto \left( \left( \begin{array}{cc} z & 0 \\ 0 & z^{-1} \\ 
\end{array} \right) U, \left( \begin{array}{cc} z^{-1} & 0 \\ 0 & z \\ 
\end{array} \right) V \right).
$$
It is clear that the orbits of this action coincide with the integral curves of the field  
$\xi_1-\tilde{\xi}_1$.  Thus, it is possible to continue the metric  $d\bar{s}^2$ on 
$[0,\infty) \times M$  by contracting each orbit into a point  for $t=0$.
The diffeomorphism 
$$
\phi: M \rightarrow M: (U, V) \mapsto (U, U^{-1} V).
$$
transforms the above-considered action of  $U(1)$ into the action of the following  form:
$$
z \in U(1): (U, V) \mapsto \left( \left( \begin{array}{cc} z & 0 \\ 0 & z^{-1} \\ 
\end{array} \right) U, V \right).
$$ 
Factorization with respect to the $U(1)$ action on the first factor defines the Hopf fibration  $G=S^3 
\rightarrow S^2=G/U(1)$. After contraction into a point of the orbits of this action for $t=0$ in the space 
$[0,\infty) \times G$, we obtain a cylinder of the Hopf fibration, which can be easily seen to be  
homeomorphic  to a linear $\mathbb{C}$-bundle  $H$ over  $S^2$, called a tautological  
bundle over $S^2$. Since the action on the second factor is trivial, we conclude that the metric
$d\bar{s}^2$ can be continued to the space $H \times G$.  

Consider now the cyclic subgroup  $\mathbb{Z}_4$  in  $U(1)$. The group  
$\mathbb{Z}_4$ (following  $U(1)$) acts on $M$; therefore, it is possible to expand this action to the 
entire space $\bar{M}$. Since this discrete action is in agreement with the orbits contraction
for $t=0$, we obtain an action of  $\mathbb{Z}_4$ on $H \times G$.  The factor-space  
$(H \times G)/\mathbb{Z}_4$ is naturally diffeomorphic  to  ${\cal 
M}=H^4 \times G$, where  $H^4$ is the fourth tensor power of the bundle 
$H$. Thus, in the considered case the metric $d\bar{s}^2$ can be continued to the manifold 
 ${\cal M}$.

The next lemma is proven analogously to the Lemma 5 from \cite{Baz}.

\vskip0.2cm

{\bf Lemma 2}. {In order for the metric $d\bar{s}^2$ to be continued to a smooth metric on ${\cal M}$, 
it is necessary and sufficient that the following conditions hold:

(1) $B_1(0)=0, |B_1'(0)|=2$;

(2) $A_2(0)=B_2(0) \neq 0, A_2'(0)=-A_2'(0)$,

(3) $A_1(0) \neq 0, A_1'(0)=0$;

(4) the functions  $A_i, B_i$  are sign-definite on the interval 
$(0,\infty)$. }

\vskip0.2cm

{\bf Remark }. A dissimilarity with the paper  \cite{Baz},  which appears  in the condition for the initial 
derivative of the function $B_1$  in (1), is connected with normalization: the length of the vector 
 $\xi_1-\tilde{\xi}_1$ equals  2, not one, as it was in Lemma 5 from \cite{Baz}.

In \cite{Brandhuber} an exact solution of the following form was found for system (2) (other solutions 
of the family, discovered in \cite{Brandhuber}, are homothetic to this one):
$$
\begin{array}{l}
A_1 (r)= \sqrt{\frac{\left(r-9/4\right) \left(r+9/4 \right) }{\left(r-3/4\right) 
\left(r+3/4 \right) }}, \\
A_2(r)= \frac{1}{\sqrt{3}}\sqrt{\left(r+3/4\right)\left(r-9/4\right)},\\
B_1(r)=2r/3, \\
B_2(r)= \frac{1}{\sqrt{3}}\sqrt{\left(r-3/4\right)\left(r+9/4\right)},
\end{array}\eqno{(3)}
$$
where  $r\geq 9/4$, and the variable  $r$ is connected with $t$ via the variables change 
$$
dt=\frac{dr}{A_1(r)}, t|_{r=\frac{9}{4}}=0.
$$
Metric (3) is a complete metric with holonomy  group $G_2$ on $S^3 
\times \mathbb{R}^4$.
If we consider the case $A_1=A_2=A_3=A$ and $B_1=B_2=B_3=B$, then the system 
(2) can be integrated in elementary functions, and we obtain another complete metric with 
holonomy group $G_2$ on $S^3 \times \mathbb{R}^4$:
$$
d\bar{s}^2=\frac{dr^2}{1-\frac{1}{r^3}}+\frac{r^2}{9} \left( 
1-\frac{1}{r^3} \right) \sum_{i=1}^3 \left( \eta_i+\tilde{\eta_i}\right)^2 + 
\frac{r^2}{3} \sum_{i=1}^3 \left( \eta_i-\tilde{\eta_i}\right)^2.\eqno{(4)}
$$
Metric (4) was constructed for the first time in \cite{Br-Sal}, see also \cite{GPP}. 
As far as we know, metrics (3) and (4) exhaust the list of known explicit solutions of the system 
(2) corresponding to complete Riemannian metrics with holonomy group $G_2$.

If we perform a formal variables change $r \rightarrow -r$ in the solution (4), then we get the following 
solution of (2):
$$
d\bar{s}^2=\frac{dr^2}{1+\frac{1}{r^3}}+\frac{r^2}{9} \left( 
1+\frac{1}{r^3} \right) \sum_{i=1}^3 \left( \eta_i+\tilde{\eta_i}\right)^2 + 
\frac{r^2}{3} \sum_{i=1}^3 \left( \eta_i-\tilde{\eta_i}\right)^2.\eqno{(5)}
$$
The solution  (5) is defined for  $0<r<\infty$, but it does not yield any smooth Riemannian metric , 
because it has a singularity at $r=0$.

\section{A family of new solutions}

Proceeding analogously to  \cite{Baz}, we consider the standard Euclidian space 
 $\mathbb{R}^4$ and set 
$R(t) = (A_1(t),  A_2(t), B_1(t), 
B_2(t)).$
Let 
 $V: \mathbb{R}^4 \rightarrow \mathbb{R}^4$  be the function of the argument   
$R$,
defined by the right-hand side of system (1) ( the function $V$ is defined, of course, only within the 
domain where  $A_i$, $B_i \neq 0$). Thus, system (1) has the form:
$$
\frac{dR}{dt}= V(R).
$$
Using the invariance of  $V$  with respect to homotheties  $\mathbb{R}^4$,
we perform substitution  $R(t)=f(t)S(t)$, where
$$
\begin{array}{c}
|S(t)|=1, f(t)=|R(t)|, \\ S(t)=(\alpha_1(t),
\alpha_2(t),\alpha_3(t),\alpha_4(t) ). \end{array}
$$
Thus, our system is split  into ''radial''  and  ''tangential'' parts:
$$
\frac{dS}{du}=V(S)-\langle V(S), S \rangle S=W(S),\eqno{(6)}
$$
$$
\begin{array}{l}
\frac{1}{f} \frac{df}{du}= \langle V(S),S \rangle, \\ dt=f
du.\end{array}\eqno{(7)}
$$
The solutions of the autonomous system (6) on the three-dimensional sphere 
$$
S^3 =\{ (\alpha_1,\alpha_2,\alpha_3, \alpha_4)|
\sum_{i=1}^4 \alpha_i^2=1\}
$$
allow us obtaining  the solutions of (2) by integrating  equations (7).
The following lemma is obvious.

\vskip0.2cm

{\bf Lemma 3}. {\it Systems (2) and (6) admit the following symmetries:
$$
\begin{array}{l}
(\alpha_1, \alpha_2, \alpha_3, \alpha_4) \mapsto
(- \alpha_1, \alpha_4, \alpha_3,\alpha_2), \\
((\alpha_1(u),\alpha_2(u), \alpha_3(u), \alpha_4(u)) \mapsto
(-\alpha_1(-u),\alpha_2(-u), \alpha_3(-u), -\alpha_4 (-u)),\\
((\alpha_1(u),\alpha_2(u), \alpha_3(u), \alpha_4(u)) \mapsto
(-\alpha_1(-u),-\alpha_2(-u), \alpha_3(-u), \alpha_4 (-u)),\\
((\alpha_1(u),\alpha_2(u), \alpha_3(u), \alpha_4(u)) \mapsto
(\alpha_1(u),\alpha_2(u), - \alpha_3(u), -\alpha_4 (u)),\\
((\alpha_1(u),\alpha_2(u), \alpha_3(u), \alpha_4(u)) \mapsto
(\alpha_1(u),-\alpha_2(u), -\alpha_3(u), \alpha_4 (u))
\end{array}
$$}

\vskip0.2cm

By virtue of Lemma 2, to the regular metric on ${\cal M}$ there may correspond only a trajectory of 
system (6) coming out of the point $S_0=(\mu, \lambda, 0, 
\lambda)$,  where $2\lambda^2+\mu^2=1$.  Due to symmetries of Lemma 3, we can assume that
 $\lambda, \mu>0$.

\vskip0.2cm

{\bf Lemma 4}. {\it For any point 
$S_0=(\mu, \lambda, 0, \lambda)$ , considered above , there exists a unique smooth trajectory 
of system (6), coming out of the point $S_0$ into the domain 
 $\alpha_3>0, \alpha_4>\alpha_2$. }

\vskip0.2cm

{\bf Proof}. Let  $J=\{ (\mu, \lambda, 0, \lambda) | \mu>0, 
\lambda>0, 2\lambda^2+\mu^2=1 \}$  be an arc of the circle on which the point $S_0$ is selected.
Let us denote by $U$ an open disc in $\mathbb{R}^2$  with coordinates 
$x=\alpha_3$, $y=\alpha_4-\alpha_2$  of the radius  $\varepsilon$  with the center at zero. Then in a 
neighborhood of the arc $J$ we can consider the local coordinates  $x, y, z=\alpha_1$. In these 
coordinates the field $W$ has the following components:
$$
W_x=W_3,\  W_y=W_2-W_4,\  W_z=W_1,
$$
where
$$
W_j(S)=V_j(S)-\langle V(S), S \rangle \alpha_j, \ j=1,2,3,4,
$$
$$
S=(\alpha_1,\alpha_2,\alpha_3,\alpha_4) =
$$
$$
\left(z, \frac{1}{2}\left(\sqrt{2-2x^2-y^2-2z^2}-y\right), x,  
\frac{1}{2}\left(\sqrt{2-2x^2-y^2-2z^2}+y\right)\right),
$$
whereas the formulae for $V_i(S)$ are obtained by the corresponding coordinate change. 
Since at the points of  $J$ the original system has a singularity, we consider in the neighborhood  
$J \times U$ a modified system of differential equations:
$$
\frac{d}{dv} \left( \begin{array}{c} x \\ y \\ z \end{array}
\right) = \left( \begin{array}{c} x W_x \\ x W_y \\ x W_z
\end{array} \right).\eqno{(5)}
$$
Clearly,  the trajectories of system (5) coincide with the trajectories of system (2) up to the parameter 
change  $du=x dv$. The vector field 
$x W$ is smooth in the neighborhood  $J \times U$;  and a direct calculation shows that for sufficiently 
small  $\varepsilon>0$  the stationary points of system (5) in  $J \times U$ are precisely 
the points of the interval $J$. Consider linearization of system (5) in a neighborhood of the point 
$S_0$:
$$
\begin{array}{l}
\frac{d x}{dv}= 2x, \\
\frac{d y}{dv}= \frac{\mu x}{\sqrt{2-2\mu^2}}-y, \\
\frac{d z}{dv}= 0.
\end{array}
$$
The linearized system has three eigenvectors 
$e_1=(3,\frac{\mu}{\sqrt{2-2\mu^2}},0)$, $e_2=(0,1,0)$, $e_3=(0,0,1)$  
with the eigenvalues $2$, $-1$ and $0$, respectively.

A direct calculation shows that if  $(x,y,z) \rightarrow 
S_0=(0,0,\mu)$, 
then
$\langle (0,0,1), \frac{x W}{|x W|} \rangle \rightarrow 0$, i.~e.
the angle between the vector $x W$ and the vector, which is tangent to the arc  $J$,
tends to $\pi/2$ when we approach the points of $J$. This allows reconstructing  the  ''phase 
portrait'' of system (5) in the neighborhood of $J \times U$ analogously to the way it is done in the 
classical case. Namely, consider the domain  $\Gamma$  in  $J \times U$, bounded by the parabolic 
cylinders 
$-\frac{\mu x}{\sqrt{2-2\mu^2}} +3y+\alpha x^2=0$, $-\frac{\mu 
x}{\sqrt{2-2\mu^2}} +3y-\alpha x^2=0$  
and the plane  $x=\delta$, where $\alpha, \delta>0$. These cylinders at the fixed level $z$
are parabolas, which are tangent along the vector $e_1$.  It is easy to calculate that at the points of the 
first parabolic cylinder 
$$
\frac{d}{dv}\left(-\frac{\mu x}{\sqrt{2-2\mu^2}} +3y+\alpha x^2\right)= 5 
\alpha x^2 + O(x^2+y^2) \geq 0,
$$
if we choose the constant  $\alpha$  to be sufficiently large (whereas equality is reached only on $J$).
Thus, the trajectories intersect the first parabolic cylinder, coming from outside of the domain 
 $\Gamma$ inside.
It can be demonstrated analogously that the trajectories of system (5) intersect the second parabolic
cylinder, bounding the domain $\Gamma$, also passing  from outside of the domain inside. Then for 
each value $z=z_0$  there exists a trajectory, which ends on the planar wall of the domain at the point  
$(\delta, y, z_0)$, and which comes out of a point on the axis $J$, if we choose $\delta$ sufficiently 
small and $\alpha$ sufficiently large (this follows from that such trajectory cannot deviate substantially 
along  $J$, since the angle that it forms with $J$ converges to $\pi/2$).  Hence, if we fix the point 
 $S_0=(0,0,\mu)$ on the arc $J$, then under diminishing of $\delta$ and increasing of  $\alpha$
we can find a trajectory, coming out exponentially with the order of $e^{2v}$  from the point  $S_0$ into 
the domain $x>0$. Analogously, there will exists a trajectory, coming out of $S_0$ from the opposite 
side, i.e. from the side of the domain $x<0$. Since the order of convergence of $x$ to zero equals 
 $e^{-2v}$, then, with respect to the parameter $u$, there will take place ''coming out'' from the 
point  $S_0$ in finite time. By analogous reasoning we demonstrate the uniqueness of each of the 
trajectories.   
Let us note now that under the transition from the parameter  $u$ to the parameter $v$
there occurs the  reversal of the trajectoriesХ orientation in the domain $x<0$. It means that for each 
point  $S_0$ there exists a unique trajectory which comes out of the point $S_0$ in finite time 
and enters the domain  $x>0$. Moreover, the trajectory, coming out of $S_0$, will be tangent to the 
vector  $e_1$, i.e. for small $u$ we will have  $\alpha_4>\alpha_2$. The Lemma is proved.

\vskip0.2cm

{\bf Lemma 5}. {\it The stationary solutions of system (6) on $S^3$ are exhausted by the following list of 
zeros of the vector field   $W$, up to the symmetries of Lemma 3:
$$
\left(\frac{1}{2\sqrt{2}}, \frac{1}{2\sqrt{2}}, \frac{\sqrt{3}}{2\sqrt{2}}, 
\frac{\sqrt{3}}{2\sqrt{2}}\right), \
\left(0, \frac{\sqrt{3}}{\sqrt{10}},\frac{\sqrt{2}}{\sqrt{5}}, 
\frac{\sqrt{3}}{\sqrt{10}} \right). \\
$$
}

\vskip0.2cm

{\bf Proof}. At the points, where the vector field $W$ turns to zero, the field
$V(S)$ is parallel to $S(u)$; hence, the stationary solutions of system (2) satisfy the following 
system of equations 
$$
\begin{array}{l}
\frac{1}{2}\left(\frac{\alpha_1^2}{\alpha_2^2}-\frac{\alpha_1^2}{\alpha_4^2}\right) 
= \beta \alpha_1,\\
\frac{1}{2}\left(\frac{\alpha_4^2-\alpha_2^2+\alpha_3^2}{\alpha_3 
\alpha_4}-\frac{\alpha_1}{\alpha_2}\right)= \beta \alpha_2,\\
\frac{\alpha_2^2+\alpha_4^2-\alpha_3^2}{\alpha_2\alpha_4}= \beta 
\alpha_3,\\
\frac{1}{2}\left(\frac{\alpha_2^2-\alpha_4^2+\alpha_3^2}{\alpha_2 
\alpha_3}+\frac{\alpha_1}{\alpha_4}\right)= \beta \alpha_4, \\
\alpha_1^2+\alpha_2^2+\alpha_3^2+\alpha_4^2=1,
\end{array}
$$
where $\beta=\langle V(S), S\rangle \in \mathbb{R}$.  Solution of the system is subdivided into two 
cases: if  $\alpha_1=0$, then we easily obtain the second point from the conditions of the Lemma. If 
 $\alpha_1 \neq 0$, then eliminating  $\beta$, we express  $\alpha_1, \alpha_3$  in terms of $\alpha_2, \alpha_4$:
$$
\alpha_1^2 = 
\frac{4}{3}\frac{\alpha_2^2\alpha_4^2}{\alpha_2^2+\alpha_4^2},
\alpha_3^2=3 
\frac{\left(\alpha_4^2-\alpha_2^2\right)^2}{\alpha_2^2+\alpha_4^2},
$$
after which, we obtain the relation 
$$
4 \left(\alpha_4^2-\alpha_2^2\right)^2 = 
\left(\alpha_4^2+\alpha_2^2\right)^2,
$$
from which we immediately obtain the remaining points . The Lemma is proved.

\vskip0.2cm

A point  $S \in S^3$,where the field  $W$ is not defined, will be called  {\it
conditionally stationary},if  there exists a real-analytic curve  
$\gamma(u)$
on $S^3$, $u \in (-\varepsilon,\varepsilon)$, $\gamma(0)=S$ , such that
the fields  $V$, $W$ are defined at all points  $\gamma(u)$, $u \in
(-\varepsilon,\varepsilon)$, $u \neq 0$, are continuously extendable to the entire curve $\gamma(u)$,
and  $\lim_{u \rightarrow 0} W(\gamma(u))=0$.

\vskip0.2cm

{\bf Lemma  6}. {\it System (6) does not have any conditionally stationary solutions on  
$S^3$.}

\vskip0.2cm 

{\bf Proof}. Let a point 
$S=(\alpha_1,\alpha_2,\alpha_3,\alpha_4)$,  $\sum_{i=1}^4
\alpha_i^2=1$ 
be conditionally stationary, i.e. there exists a curve 
$\gamma(u)$, $u \in (-\varepsilon,\varepsilon)$ 
with the above-mentioned properties. Obviously, this is only possible in the case when at least one of the conditions holds:
 $\alpha_2(0)=0$, $\alpha_3(0)=0$  or 
$\alpha_4(0)=0$.

1) First consider the case when all the relations hold simultaneously:
 $\alpha_2(0)=\alpha_3(0)=\alpha_4(0)=0$, $\alpha_1(0)=\pm 
1$. 
Let us set for $i=2, 3, 4$
$$
\alpha_i(u)=c_i u^{k_i} (1+o(1)), \mbox{\ при\ } u \rightarrow 0,
$$
where
 $c_i \neq 0$, $k_i >0$.  Note that if 
 $\alpha_2(u)=\alpha_4(u)+c 
u^{k}$,  where $c \neq 0$, then  $V_1$ cannot be continuously extended along  
 $\gamma(u)$ up to $u=0$. It follows from the real analyticity that 
$\alpha_2(u)=\alpha_4(u)$ and, in particular, $k_2=k_4$.  Then 
$$
V_2=\frac{\pm 1}{2 c_2}u^{-k_2}(1+o(1)),
$$
which is a contradiction with the existence of the limit of $V(S)$ as $u \rightarrow 0$.

2) Suppose that two out of three functions  $\alpha_2, \alpha_3, \alpha_4$  turn to zero at $u=0$. 
Consider the arising cases.

The case of $\alpha_2(0)=\alpha_3(0)=0$, $\alpha_4(0)\neq 0$. If, in addition, 
$\alpha_1(0)\neq 0$,  then
$$
V_1=\frac{\alpha_1(0)}{2c_2^2}u^{-2k_2} (1+o(1)),
$$
which leads to a contradiction. If  $\alpha_1(u)=c_1 u^{k_1} (1+ o(1))$, 
$c_1\neq 0$, $k_1>0$, then $k_1 \geq k_2$ (from the continuity of $V_1$) and
$$
V_2= \frac{\alpha_4(0)}{c_3}u^{-k_3} (1+o(1)),
$$
which is again a contradiction.

The case of $\alpha_2(0)\neq 0$, $\alpha_3(0)=\alpha_4(0)=0$ is symmetric to the previous one and 
can be excluded analogously.  

The case of $\alpha_2(0)=\alpha_4(0)=0$, $\alpha_3(0)\neq 0$. This case is excluded, because  
$$
V_3=-\frac{\alpha_3(0)^2}{c_2 c_4} u^{-k_2-k_4} (1+o(1)).
$$

3) Suppose that only one of the functions $\alpha_2, \alpha_3, \alpha_4$ 
turns to zero at $u=0$.

The case of $\alpha_2(0)=0$, $\alpha_3(0), \alpha_4(0) \neq 0$. The continuity of  
$V_1$ and $V_3$ implies in this case that  $\alpha_1(0)=0$ and $\alpha_3(0)=\pm 
\alpha_4(0) \neq 0$, while $k_1\geq k_2$. In this case  $\lim_{u\rightarrow 
0} V(\gamma(u))=(0, 1, 0, 0)$ and  $\lim_{u\rightarrow 0} W(\gamma(u))=(0, 1, 
0, 0) \neq 0$, which is a contradiction.

The case of  $\alpha_4(0)=0$, $\alpha_2(0), \alpha_3(0) \neq 0$ is excluded in analogous fashion.

The case of  $\alpha_3(0)=0$, $\alpha_2(0), \alpha_4(0) \neq 0$. The continuity of  
$V$ immediately yields that $\alpha_2(0)=\pm \alpha_4(0)$. Then  $\lim_{u\rightarrow 0} 
V_3(\gamma(u))=2$  and  $\lim_{u\rightarrow 0} W_3(\gamma(u))=2 \neq 0$, again a contradiction.
 The Lemma is proved.

A metric  $d\bar{s}^2$ is called asymptotically locally conic, if there exist the functions 
$\tilde{A}_i (t)$, $\tilde{B}_i(t)$, linear with respect to $t$ up to a shift,  such that 
$$
\left|1-\frac{A_i}{\tilde{A}_i} \right| \rightarrow 0, \ 
\left|1-\frac{B_i}{\tilde{B}_i} \right| \rightarrow 0, \mbox{\ при \ } t 
\rightarrow \infty
$$
The metric, defined by the functions $\tilde{A}_i (t)$, $\tilde{B}_i(t)$,  is called locally conic. The following lemma is proved in  \cite{Baz}. 

\vskip0.2cm

{\bf Lemma 7}. {\it To the stationary solutions of system (6) there correspond the locally conic metrics 
on $\overline{M}$, whereas to the trajectories of system (6), asymptotically converging to stationary 
solutions,  there correspond the asymptotically locally conic metrics on 
$\overline{M}$. 
}

\vskip0.2cm

The following lemma follows directly from the analysis of systems (2) and (6).

\vskip0.2cm

{\bf Lemma 8}. {\it  If $S=(\alpha_1,\alpha_2,\alpha_3,\alpha_4)$
is a solution of system (6), then there take place the following relations:
$$
\begin{array}{ll}
1) & \frac{d}{dt} \left( 2 A_1 A_2 B_2 - B_1 ( B_2^2-A_2^2)\right) =0, \\
2) & \frac{d}{du} \left( \frac{\alpha_1 \alpha_2 
\alpha_4}{2\alpha_4\alpha_2\alpha_1-\alpha_3 \left( \alpha_4^2-\alpha_2^2 
\right)}\right) = \frac{\alpha_1\alpha_3}{2\alpha_4\alpha_2\alpha_1-\alpha_3 
\left( \alpha_4^2-\alpha_2^2 \right)}, \\
3) & \frac{d}{du} \left( \ln \frac{\alpha_3 \left( \alpha_4^2-\alpha_2^2 
\right)}{\alpha_4\alpha_2\alpha_1} \right)=
\frac{2\alpha_4\alpha_2\alpha_1-\alpha_3 \left( \alpha_4^2-\alpha_2^2 
\right)}{2\alpha_4\alpha_2\left( \alpha_4^2-\alpha_2^2 \right)}, \\
4) & \frac{d}{du} \ln \frac{\alpha_2}{\alpha_4} = 
\frac{\alpha_4^2-\alpha_2^2}{\alpha_2\alpha_3\alpha_4}, \mbox{\ при \ } 
\alpha_2=\alpha_4, \\
5) & \frac{d}{du} \left( \frac{\alpha_3}{\alpha_4}\right) = 
\frac{3}{2\alpha_4} 
\left(\frac{2}{\sqrt{3}}+\frac{\alpha_3}{\alpha_4}\right)\left(\frac{2}{\sqrt{3}} 
-\frac{\alpha_3}{\alpha_4}\right) \mbox{\ при \ } \alpha_1=0, 
\alpha_2=\alpha_4.
\end{array}
$$
}

\vskip0.2cm

{\bf Remark}.  Thus, the function  
$F(t)=2A_1A_2B_2-B_1(B_2^2-A_2^2)$  is an integral of system (2).

\vskip0.2cm

{\bf Lemma 9}. {\it The trajectory of system (6), defined by the initial point  
$S_0=(\mu,\lambda,0,\lambda)$, $\lambda,\mu >0$, $2\lambda^2+\mu^2=1$, 
converges, as $u\rightarrow \infty$, to the stationary point 
$S_\infty=\left(0, \frac{\sqrt{3}}{\sqrt{10}}, \frac{\sqrt{2}}{\sqrt{5}}, 
\frac{\sqrt{3}}{\sqrt{10}}\right)$}.

\vskip0.2cm

{\bf Proof}. Let us introduce notations for the following points in $S^3$:
$$
\begin{array}{l}
O=(0,0,1,0), \
A=(0,0,0,1), \
B=(1,0,0,0), \
C=(0, \frac{1}{\sqrt{2}}, 0, \frac{1}{\sqrt{2}}).
\end{array}
$$
Consider the domain $\Pi \subset S^3$, defined by the inequalities:
$$
\Pi: \alpha_4\geq \alpha_2\geq0, \alpha_1 \geq0, \alpha_3 \geq0.
$$
It is not difficult to verify that the domain  $\Pi$ is the spherical pyramid  
$(OABC)$.
The boundaries of the domain are the following sets:
$$
\begin{array}{l}
\Pi_1=(OAB)=\{ \alpha_2=0, \alpha_4\geq 0, \alpha_1\geq0, \alpha_3\geq0\} , 
\\
\Pi_2=(OBC)=\{ \alpha_4= \alpha_2, \alpha_2\geq0, \alpha_1\geq0, 
\alpha_3\geq0 \}, \\
\Pi_3=(OAC)=\{ \alpha_4 \geq \alpha_2 \geq 0, \alpha_1=0, \alpha_3\geq0\}, 
\\
\Pi_4=(ABC)=\{ \alpha_4 \geq \alpha_2 \geq0, \alpha_1 \geq0, \alpha_3=0 \}.
\end{array}
$$
The initial point  $S_0=(\mu, \lambda, 0, \lambda) \in (BC)$. According to Lemma 4, for all small 
 $u>0$  the trajectory of system (6), determined by the initial point  $S_0$, is inside the domain 
 $\Pi$.

Consider first the possibility of the trajectory reaching the boundary of the domain  
$\Pi$ in finite time. On $\Pi_1\backslash ((AB)\cup (OB))$ the integral  
$F(t)=-\alpha_3\alpha_4^2 f(t)^3$  is strictly negative, whereas at the initial point  
$F(S_0)=2\lambda^2\mu>0$; hence, the trajectory cannot intersect a certain neighborhood of this wall,
with the possible  exception of the arcs  $(AB)$ and  $(OB)$. 
Further, on $\Pi_2$  we have 
$$
\frac{d (\alpha_4-\alpha_2)}{du} = \frac{\alpha_1}{\alpha_2}>0,
$$
for  $\alpha_1 \neq 0$, i.e.  the trajectory cannot interest a certain neighborhood of the set 
$\Pi_2$ in finite time, or even come sufficiently close to it, with the exception of the arc 
  $(OC)$. Notice that this consideration also excludes a neighborhood of the arc $(OB)$. Finally, on the 
set  $\Pi_4$  the derivative of the function  $\alpha_3(t)$  is strictly positive and bounded away from
 zero, therefore the trajectory does not intersect  $\Pi_4$ and its certain neighborhood  (note that  we
 have excluded at once the remaining possibility of approaching  the arc  $(AB)$). Since 
 $\Pi_3$ is an invariant subset of system (6), then the trajectory cannot intersect  
$\Pi_3$ in finite time (including the arc $(OC)$).

Define a function $F_1$ on $S^3$:
$F_1(\alpha_1,\alpha_2,\alpha_3,\alpha_4)=
\frac{\alpha_1 \alpha_2 \alpha_4 }{F(\alpha_1, \alpha_2, \alpha_3, 
\alpha_4)}$.
Since $F(\alpha_1,\alpha_2,\alpha_3,\alpha_4)=f(t)^{-3} F(S_0)>0$, then it follows from relation 2)
of  Lemma 8 that the function $F_1$ strictly increases along the trajectories of system (2), passing 
inside the domain $\Pi$. Suppose that 
$C$  is the limit set of the considered trajectory. Then, only the following points can get into $C$: either 
stationary and conditionally stationary points of system (6)  
(i. e., according to Lemmas  5 and 6 we have only two such possibilities: the points  $S_\infty$ and 
 $S_1=\left(\frac{1}{2\sqrt{2}}, 
\frac{1}{2\sqrt{2}}, \frac{\sqrt{3}}{2\sqrt{2}},  
\frac{\sqrt{3}}{2\sqrt{2}}\right)$);  or the points lying on the critical level of the function $F_1$ (it is clear 
that there are no such points inside $\Pi$, since in the neighborhood of each point of $C$ interior with 
respect to $\Pi$ it is possible to bound away from zero the derivate of 
$F_1(u)$); finally, all the points of $C$, lying on the boundary of 
$\Pi$, must be situated on the maximal level of the function $F_1$. Quite analogously, we consider 
the function  
$F_2(\alpha_1,\alpha_2,\alpha_3,\alpha_4)=
\ln \frac{\alpha_3 \left( \alpha_4^2-\alpha_2^2 
\right)}{\alpha_4\alpha_2\alpha_1}$.  Then it follows from the relation 3) of Lemma 8 that  
$F_2$ is increasing along the trajectory; and thus, the set $C \cap 
\partial \Pi$ lies on the maximal level of  $F_2$ in $\Pi$.  Let us note that the maximal  
(in $\Pi$) level of the function $F_2$ is the set  $\Pi_3 
\cup \Pi_1$. It has been demonstrated above that it is impossible to approach a neighborhood of 
$\Pi_1 \backslash (OA)$; hence, the case of $C \cap 
\partial \Pi \subset \Pi_3$ is the only possible one.

Now, it follows from the relation  4) of Lemma 8 that the function $F_3=\ln 
\frac{\alpha_2}{\alpha_4}$  is increasing along the trajectory  (for sufficiently large $u$) towards
the maximal value on $\Pi_3$, which is reached for 
$\alpha_2=\alpha_4$. Thus, our trajectory is converging, as $u \rightarrow 
\infty$, to the invariant one-dimensional set $\Pi_3\cap \Pi_2=(OC)$.
The relation 5) of Lemma 8 shows that, in the neighborhood of $(OC)$, the function  
$F_4=\frac{\alpha_3}{\alpha_4}$ is increasing for  $F_4 \leq 
\frac{2}{\sqrt{3}}$ and decreasing for $F_4 \geq \frac{2}{\sqrt{3}}$; 
hence,  $C \cap \partial \Pi$ can contain only the point  $S_\infty$, 
defined by the condition $F_4 =\frac{2}{\sqrt{3}}$.

Thus, we have arrived at the conclusion that the considered trajectory converges either to  
$S_1$, or to $S_\infty$. To complete the proof, it remains  for us to show that  the convergence to 
$S_1$ does not take place.

A direct calculation shows that the linearization of system (6) in the neighborhood of the stationary 
point  $S_1$ has three eigenvalues of multiplicity one:
$$
\lambda_1=-2 \sqrt{2}, \lambda_2 = -\frac{7}{3}\sqrt{2} - 
\frac{1}{3}\sqrt{290}, \lambda_3= -\frac{7}{3}\sqrt{2} + 
\frac{1}{3}\sqrt{290}.
$$
Thus, in the neighborhood of the point  $S_1$ there exists a (locally defined ) surface, formed
by the trajectories, entering the point  
$S_1$; moreover, this surface is tangent at the point $S_1$ to the two-dimensional plane
spanned by the first two eigenvectors  $e_1$ and $e_2$. Other trajectories in the neighborhood of  
 $S_1$ come out of $S_1$. At that, the first eigenvector has the following coordinates (in  
$\mathbb{R}^4$):  $e_1=(-\sqrt{3},-\sqrt{3},1,1)$ and  is tangent to the trajectory which is defined as 
 $\alpha_1=\alpha_2$, $\alpha_3=\alpha_4$. It is not difficult to see that to the eigenvalue  
$\lambda_1$ there correspond precisely the solutions (4) and (5) 
(these trajectories enter the point  $S_1$ from the opposite sides; trajectory  
(4) corresponds to  $F<0$, whereas trajectory (5), to $F>0$). Since  
$|\lambda_2|>|\lambda_1|$, then other trajectories entering $S_1$  
(with the exception of one of them) are tangent at the point $S_1$ to the trajectory (5) or (6). The only
non-tangent  to (5), (6)  trajectory, mentioned above, corresponds to the eigenvalue      
$\lambda_2$; and it can be directly checked that it lies on the invariant surface 
$F=0$ and, thus, cannot coincide with our trajectory.

Consider a couple of functions: $G_1=\alpha_2\alpha_4-\alpha_1\alpha_3$ and  
$G_2=\alpha_1\alpha_4-\alpha_2\alpha_3$. The initial point  $S_0$ is situated in the region  
$\{G_1>0, G_2>0\}$, the point $S_1$ lies in  $\{G_1=0, G_2=0\}$. 
A direct calculation shows that the vector $e_2$ is directed inside the domains $\{G_1>0, G_2>0\}$ or 
 $\{G_1<0, G_2<0\}$ (the domain depends on the choice of direction of  $e_2$; to avoid lengthy 
formulae, we do not present here the explicit coordinates of $e_2$ ). It is easy to check that   
$\frac{d}{du}G_1=-\frac{2}{\alpha_2} G_2$,  at the points where $G_1=0$;  and 
$\frac{d}{du}G_2=- \frac{2}{\alpha_2} G_1$ at those points where $G_2=0$. 
Hence, the trajectory can reach the point $S_1$ only by staying in the domain  
$\{G_1>0, G_2>0\}$; if it moves to one of the regions $\{G_1>0, 
G_2<0\}$ or $\{G_1<0, G_2>0\}$, then it will not be able to leave them  
(let us note that $S_\infty$ lies in $\{ G_1>0, G_2<0\}$). Among other things, this consideration 
determines the direction of the vector $e_2$: it is directed inside the domain  $\{G_1>0, G_2>0\}$.

Consider now the function $F_5=\alpha_4^2-\alpha_3^2$. It is obvious that  
$F_5(S_0)=\lambda^2>0$, $F_5(S_1)=0$.  Next,
$$
\frac{d}{du} F_5= \frac{G_2}{2\alpha_1\alpha_4},
$$
at those points where $F_5=0$. Thus, on the level surface  
$\{F_5=0\}$  the derivative of the function along the trajectory is nonnegative and 
turns to zero precisely at the points where $\alpha_1=\alpha_2$ and 
$\alpha_3=\alpha_4$. Since these points belong to the trajectory of solution  
(5), the considered trajectory cannot leave the domain  
$\{F_5>0\}$. On the other hand, a direct computation shows that the vector $e_2$ (along which the 
trajectories come to the point $S_1$) is directed inside the domain  $\{F_5<0\}$. Thus, while staying  
in the domain $\{G_1>0, 
G_2>0\}$,  the trajectory cannot come close to $S_1$. There remains just one possibility: 
going out to the domain  $\{ G_1>0, G_2<0\}$, where the only limit point is 
$S_\infty$.
The Lemma is proved.

The main theorem is now a direct consequence of Lemmas 4 and 9: 
the initial point of the trajectory determines the topological structure of the space, on which our metric 
is defined whose holonomy group, obviously, coincides with the entire   
$G_2$. The limit point  $S_\infty$ means that the function $B_1$ is approximated at infinity by a 
constant , whereas other functions, defining the metrics, by some linear non-constant functions. At 
infinity this yields  the product of  $S^1$ and a cone over $S^2 \times 
S^3$.

{\bf Acknowledgement}. The work was supported by the Russian Foundation for Basic Research (grants no. 12-01-00873, 12-01-92104-YaF-a (first author) and 12-01-00124 (second author)), 
the Grant Council of the President of the Russian Federation  (grants  MD-249.2011.1 and NSh-544.2012.1).


\end{document}